\documentclass{article}
\usepackage{amsmath, amsfonts, amssymb, amsthm}
\newtheorem{theorem}{Theorem}
\newtheorem{prob}{Problem}

\newtheorem*{rem}{Remark}
\theoremstyle{definition}
\newtheorem*{ex}{Example}


\begin{document}

\author{Yaroslav Abramov
\footnote{The author is partially supported by
AG Laboratory SU-HSE, RF government 
grant, ag. 11.G34.31.0023}
\footnote{Laboratory of Algebraic Geometry, SU-HSE, 
7 Vavilova Str., Moscow, Russia, 117312} 
\footnote{$zroslav@gmail.com$}}

\title{Systems of resultants}
\date{}

\maketitle

\begin{abstract}

Writing down convenient explicit formulas for systems of resultants is an important but essentially 
open problem. In this paper I'll give such a formula derived from the ordinary multivariate resultant. 

\end{abstract}

\paragraph*{Acknowledgements}

I'm very thankful for A.L. Gorodentsev
for useful discussions. The author is partially supported by
AG Laboratory SU-HSE, RF government 
grant, ag. 11.G34.31.0023

\section{Intro}

Fix some algebraically closed field $\Bbbk$.

\begin{prob}
Given a system of polynomial equations
\begin{equation}\left\{\begin{array}{l}
f_0(x_0, \ldots , x_n)=0\\
\ldots \\
f_m(x_0, \ldots , x_n)=0
\end{array}\right.\label{sys1}\end{equation}
$\deg f_i=N_i$, 
$$f_j(x_0, \ldots , x_n)=\sum\limits_{\sum\limits_i s_i=N_j}a_{j,s_0, \ldots , s_n}
x_0^{s_0}\cdot\ldots\cdot x_n^{s_n}.$$
How to determine if there exists a non-zero solution of (\ref{sys1})?
\end{prob}

It is well-known after \cite{wanderwarden} that 
there exists a finite set of polynonials on $a_{j,s_0,\ldots , s_n}$ with integer coefficients 
$R_l(a)\in\mathbb Z[a_{j,s_0, \ldots , s_n}]_{j, s_0, \ldots , s_n}$, such that
$$
\mbox{(there exists a non-zero solution of (\ref{sys1}))} \Longleftrightarrow \forall l\ R_l(a)=0
$$

Such a set of polynomials ($R_l(a)$) is called a {\bfseries system of resultants}.

\begin{ex}
Let $\deg f_j=1$, $j=0, \ldots , m$, $f_j(x)=\sum_i a_{ji}x_i$. Then the system of resultants is
the set of maximal minors of matrix 
$$
\begin{pmatrix}
a_{00} & \ldots & a_{0n} \\
\ldots & \ldots & \ldots \\
a_{m0} & \ldots & a_{mn}
\end{pmatrix}
$$
\end{ex}

\begin{prob}
Given a system of polynomial equations
\begin{equation}\left\{\begin{array}{l}
f_0(x_0, \ldots , x_n)=0\\
\ldots \\
f_n(x_0, \ldots , x_n)=0
\end{array}\right.\label{sys2}\end{equation}
$$\deg f_i=N_i,$$ 
$$f_j(x_0, \ldots , x_n)=\sum\limits_{\sum\limits_i s_i=N_j}a_{j,s_0, \ldots , 
s_n}x_0^{s_0}\cdot\ldots\cdot x_n^{s_n}.$$
How to determine if there exists a non-zero solution of (\ref{sys2})?
\end{prob}

It is also well-known (see  \cite{gzk} for a modern explanation) that
there exist an irreducible polynonial on $a_{j,s_0,\ldots , s_n}$ with integer coefficients 
$$R(a)\in\mathbb Z[a_{j,s_0, \ldots , s_n}]_{j, s_0, \ldots , s_n},$$ such that
$$
\mbox{(there exists a non-zero solution of (\ref{sys2}))} \Longleftrightarrow R(a)=0
$$

Such a polynomial ($R(a)$) is called a {\bfseries resultant} and also denoted as $R(f_0, \ldots , f_n)$

\begin{ex}
Let $f_i(x)=\sum_j a_{ij}x_j$ then $R(f_0, \ldots , f_n)=\det(a_{ij})$.
\end{ex}

\begin{ex}
Let $f(x,y)=\sum_i a_ix^iy^{n-i}$, $g(x,y)=\sum_i b_ix^iy^{m-i}$, $a_0\ne 0$, $b_0\ne 0$. 
Then $$R(f(x,y),g(x,y))=Res(f(z,1), g(z,1))$$ where $Res$ is a famous Sylvester determinant.
$$\det\underbrace{\left(\!\!\begin{array}{ccccccc}
 a_0&a_1&\dots&a_n\\
 &a_0&a_1&\dots&a_n\\
 &&\ddots&\ddots&&\ddots\\
 &&&a_0&a_1&\dots&a_n\\
 b_0&b_1&\dots&b_m\\
 &b_0&b_1&\dots&b_m\\
 &&\ddots&\ddots&&\ddots\\
 &&&b_0&b_1&\dots&b_m
 \end{array}\!\!\right)
 }_{m+n}\!\!\!
 \begin{array}{c}
 \left.\vphantom{\begin{array}{ccccccc}
 a_0&a_1&\dots&a_n\\
 &a_0&a_1&\dots&a_n\\
 &&\ddots&\ddots&&\ddots\\
 &&&a_0&a_1&\dots&a_n
 \end{array}
 }\right\}{\scriptstyle m}
 \\
 \left.\vphantom{\begin{array}{ccccccc}
 b_0&b_1&\dots&b_m\\
 &b_0&b_1&\dots&b_m\\
 &&\ddots&\ddots&&\ddots\\
 &&&b_0&b_1&\dots&b_m
 \end{array}
 }\right\}{\scriptstyle n}
 \end{array}
 $$
\end{ex}

\section{Results on resultants}

Consider the system

\begin{equation}\left\{\begin{array}{l}
f_0(x_0, \ldots , x_n)=0\\
\ldots \\
f_m(x_0, \ldots , x_n)=0
\end{array}\right.\label{sys3}\end{equation}

$$\deg f_j=n_j,$$ 
$$f_j(x_0, \ldots , x_n)=\sum\limits_{\sum\limits_i s_i=n_j}a_{j, s_0, \ldots , s_n}.$$

Fix some positive integer numbers 
$m_i,\ i=0, \ldots , n,$ 
and 
$k_{ij},\ i=0, \ldots , n; j=0, \ldots , m,$ 
such that $m_i=k_{ij}+n_j.$ 
Consider polynomials 
$$A_{ij}(x_0, \ldots , x_n)=
\sum_{\sum\limits_l s_l = k_{ij}} b_{i, j, s_0, \ldots , s_n}x_0^{s_0}\cdot\ldots\cdot x_n^{s_n}$$ 
with 
indeterminate coefficients 
$b_{i, j, s_0, \ldots , s_n}.$

I will consider 
$$R(\sum_{j=0}^m A_{0j} f_j, \sum_{j=0}^m A_{1j}f_j, \ldots , \sum_{j=0}^m A_{nj}f_j)$$ 
as a polynomial in 
$b_{i,j, s_0, \ldots , s_n}$ for various $i, j, s_0, \ldots , s_n$.

\begin{theorem}
System (\ref{sys3}) has a non-zero solution iff
$$R=R(\sum_{j=0}^m A_{0j} f_j, \sum_{j=0}^m A_{1j}f_j, \ldots , \sum_{j=0}^m A_{nj}f_j)\equiv 0$$
as a polynomial in the coefficients 
$b_{i,j, s_0, \ldots , s_n}$ 
of $A_{ij}$. 
Thus, coefficients of $R$ form the system of resultants of $f_0, \ldots , f_m$.
\end{theorem}

\begin{ex}
Let $f_j(x)=\sum_{i}a_{ij}x_i$ and $\deg A_{ij}=0$ then $$R(\sum_{j=0}^m A_{ij}f_j)_{i=0}^n=
\sum_{J\subset \{0, \ldots m\}, \mid J\mid=n+1}\det (a_{ij})_{i=0, \ldots , n, j\in J}\prod_{j\in J}b_j.$$
\end{ex}

\paragraph*{Proof.} Assume the contrary. For $x\in\mathbb P^n$ I put 
$$H_x=\{(A_{ij})_{ij}\mid \sum_{j=0}^m A_{ij}(x)f_j(x)=0, i=0, \ldots , n\}.$$
The condition $R\equiv 0$ is equivalent to 
$$\bigoplus_{i=0}^n\bigoplus_{j=0}^m S^{k_{ij}}(\Bbbk^{n+1})=\bigcup_{x\in\mathbb P^n} H_x$$
If $x$ is not a solution of $(\ref{sys3})$ then $H_x$ is a codimension $(n+1)$ 
linear subspace in 
$$V=\bigoplus_{i=0}^n\bigoplus_{j=0}^m S^{k_{ij}}(\Bbbk^{n+1}).$$
If there are no non-zero solutions of (\ref{sys3}) then $V$ is a union
of $n$-parametric family of codimension $(n+1)$ subspaces. We get the contradiction. 

\begin{rem}
In \cite{gzk} there is a definition of mixed resultant for sections of very ample
linear bundles $L_0, \ldots , L_n$ on a dimension $n$ projective variety. 
Theorem 1 can be generalised to the case of sections of very ample linear bundles 
$$f_j\in H^0(X, L_j), j=0 \ldots , m$$ 
on a dimension $n$ projective variety $X$. Consider a system of very ample line bundles 
$C_{ij}$, $0\leq i\leq n$, $0\leq j \leq n$, s.t. $B_i=C_{ij}\otimes L_j$ for all $i, j$. 
Then the system of resultants is just the collection of coefficients
of $$R(\sum_{j=0}^m A_{ij}\otimes f_j)_{i=0}^n$$
considered as a polynomial in indeterminate 
$$A_{ij}\in H^0(X, C_{ij}).$$
\end{rem}

\begin{rem}
We get only the set-theoretical (not the scheme-theoretical) system of resultants.
\end{rem}

There are also some related results (which may be used for simplification of calculations
and which can be proved by almost exactly the same prooftext):

\begin{theorem}
Let $\deg f_0\geq \deg f_1\geq\ldots\geq\deg f_m$ and $k_{ij}=\deg f_i-\deg f_j$. Then 
system (\ref{sys3}) has a non-zero solution iff
$$R(f_0+\sum_{j=n+1}^m A_{0j}f_j, f_1+\sum_{j=n+1}^m 
A_{1j}f_j, \ldots , f_n+\sum_{j=n+1}^m A_{nj}f_j)\equiv 0$$
as a polynomial on coefficients of $A_{ij}$.
\end{theorem}

Consider vector subspaces $V_{ij}$ of 
$S^{k_{ij}}(\Bbbk^{n+1})$,
such that  
$$\{A_{ij}(x)\mid A_{ij}\in V_{ij}\}=\Bbbk$$
for all $x\in\Bbbk^{n+1}$.

\begin{ex}
$V_{ij}=S^{k_{ij}}(\Bbbk^{n+1})$
\end{ex}

\begin{ex}
$V_{ij}=\{\sum\limits_{l=0}^na_{lij}x_l^{k_{ij}}\}$
\end{ex}

\begin{theorem}
System (\ref{sys3}) has a non-zero solution iff
$$R=R(\sum_{j=0}^m A_{0j} f_j, \sum_{j=0}^m A_{1j}f_j, \ldots , \sum_{j=0}^m A_{nj}f_j)\equiv 0,$$
(where $A_{ij}\in V_{ij}$)
as a polynomial on $\bigoplus_{i=0}^n V_i$. 
Thus, coefficients of $R$ form the system of resultants of $f_0, \ldots , f_m$.
\end{theorem}

\begin{rem}
Theorem 3 is a generalisation of Theorem 1.
\end{rem}

Consider vector subspaces
$V_i$ of $\bigoplus_{j=0}^m S^{k_{ij}}(\Bbbk^{n+1})$,
such that 
$$\{(A_{i0}(x), A_{i1}(x), \ldots , A_{im}(x)\mid 
(A_{0m}, A_{1m}, \ldots , A_{im})\in V_i\}=\Bbbk^{m+1}$$
for all $x\in\Bbbk^{n+1}$ 

\begin{ex}
$V_i=\bigoplus_{j=0}^m S^{k_{ij}}(\Bbbk^{n+1})$
\end{ex}

\begin{ex}
$V_i=\bigoplus_{j=0}^m V_{ij}$
\end{ex}

\begin{ex}
$V_i=\{(\sum\limits_{l\ne 0}a_{li0}x_l^{k_{i0}}+bx_0^{k_{i0}}, 
\ldots \sum\limits_{l\ne n}a_{lin}x_l^{k_{in}}+bx_n^{k_{in}}, 
\sum\limits_{l=0}^na_{li(n+1)}x_l^{k_{i(n+1)}}, 
\ldots , \sum\limits_{l=0}^na_{lim}x_l^{k_{im}})\}$
\end{ex}

\begin{theorem}
System (\ref{sys3}) has a non-zero solution iff
$$R=R(\sum_{j=0}^m A_{0j} f_j, \sum_{j=0}^m A_{1j}f_j, \ldots , \sum_{j=0}^m A_{nj}f_j)\equiv 0$$
(where $(A_{i0}, A_{i1}, \ldots , A_{im})\in V_i$)
as a polynomial on $\oplus_{i=0}^n V_i$. 
Thus, coefficients of $R$ form the system of resultants of $f_0, \ldots , f_m$.
\end{theorem}

\begin{rem}
Theorem 4 is a generalisation of Theorem 3.
\end{rem}


\end{document}